\documentclass[amssymb, 11pt]{amsart}
\usepackage{latexsym}

\newlength{\standardunitlength}
\newtheorem{prop}{Proposition}[section]

\newtheorem{lemma}[prop]{Lemma}
\newtheorem{cor}[prop]{Corollary}
\newtheorem{theorem}[prop]{Theorem}

\begin{document}

\begin{center}
{\bf Stein's Method and Plancherel Measure of the Symmetric Group}
\end{center}

\begin{center}
{\bf Running head: Stein's Method and Plancherel Measure}
\end{center}

\begin{center}
By Jason Fulman
\end{center}

\begin{center}
University of Pittsburgh
\end{center}

\begin{center}
Department of Mathematics
\end{center}

\begin{center}
301 Thackeray Hall
\end{center}

\begin{center}
Pittsburgh, PA 15260
\end{center}

\begin{center}
Email: fulman@math.pitt.edu
\end{center}

{\bf Abstract}: We initiate a Stein's method approach to the study of the Plancherel measure of the symmetric group. A new proof of Kerov's central limit theorem for character ratios of random representations of the symmetric group on transpositions is obtained; the proof gives an error term. The construction of an exchangeable pair needed for applying Stein's method arises from the theory of harmonic functions on Bratelli diagrams. We also find the spectrum of the Markov chain on partitions underlying the construction of the exchangeable pair. This yields an intriguing method for studying the asymptotic decomposition of tensor powers of some representations of the symmetric group.

\begin{center}
To appear in Trans. AMS. Submitted 5/28/03, Minor Revisions: 7/7/03, 11/8/03
\end{center}

\begin{center}
2000 Mathematics Subject Classification: 05E10 (primary), 60C05 (secondary)
\end{center}

\begin{center}
Key words and phrases: Plancherel measure, Stein's method, character
ratio, Markov chain. \end{center}

\newpage

\section{Introduction}

	It is accurate to state that one of the most important developments in probability theory in the last century was Stein's method. Stein's method is a highly original technique and has been useful in proving normal and Poisson approximation theorems in probability problems with limited information such as the knowledge of only a few moments of the random variable. Stein's method can be difficult to work with and often the bounds arising are not sharp, even in simple problems. But it sometimes is the only option available and has been a smashing success in Poisson approximation in computational biology. Good surveys of Stein's method (two of them books) are \cite{AGG},\cite{BHJ},\cite{St1},\cite{St2}.

	Next let us recall the Plancherel measure of the symmetric
	group. This is a probability measure on the irreducible
	representations of the symmetric group which chooses a
	representation with probability proportional to the square of
	its dimension. Equivalently, the irreducible representations
	of the symmetric group are parameterized by partitions
	$\lambda$ of $n$, and the Plancherel measure chooses a
	partition $\lambda$ with probability $\frac{n!}{\prod_{x \in
	\lambda} h(x)^2}$ where the product is over boxes in the
	partition and $h(x)$ is the hooklength of a box. The
	hooklength of a box $x$ is defined as 1 + number of boxes in
	same row as x and to right of x + number of boxes in same
	column of x and below x. For example we have filled in each
	box in the partition of 7 below with its hooklength
\[ \begin{array}{c c c c}
                \framebox{6}& \framebox{4}& \framebox{2}& \framebox{1}
                \\ \framebox{3}& \framebox{1}&& \\ \framebox{1} &&&
                \end{array}, \] and the Plancherel measure would
                choose this partition with probability
                $\frac{7!}{(6*4*3*2)^2}$. Recently there has been
                significant interest in the statistical properties of
                partitions chosen from Plancherel measure. As it is
                beyond the scope of this paper to survey the topic, we
                refer the reader to the surveys \cite{AD}, \cite{De}
                and the seminal papers \cite{J}, \cite{O}, \cite{BOO}
                for a glimpse of the remarkable recent work on
                Plancherel measure.
	
	A purpose of the present paper is to begin the study of Plancherel measure by Stein's method. A long term goal of this program is to use Stein's method to understand the Baik-Deift-Johansson theorem, giving explicit bounds on the convergence of the first row of a Plancherel distributed partition to the Tracy-Widom distribution. We can not at present do this but are confident that the exchangeable pair in this paper is the right one. We do attain a more modest goal of a Stein's method approach to the following result of Kerov. 

\begin{theorem} \label{ker} (\cite{K1}) Let $\lambda$ be a partition of $n$ chosen
 from the Plancherel measure of the symmetric group $S_n$. Let
 $\chi^{\lambda}(12)$ be the irreducible character of the symmetric
 group parameterized by $\lambda$ evaluated on the transposition
 $(12)$. Let $dim(\lambda)$ be the dimension of the irreducible
 representation parameterized by $\lambda$. Then the random variable
 $\frac{n-1}{\sqrt{2}} \frac{\chi^{\lambda}(12)}{dim(\lambda)}$ is
 asymptotically normal with mean 0 and variance 1. \end{theorem}

	Let us make some remarks about Theorem \ref{ker}. The quantity
	$\frac{\chi^{\lambda}(12)}{dim(\lambda)}$ is called a
	character ratio and is crucial for analyzing the random walk
	on the symmetric group generated by transpositions
	\cite{DS}. In fact Diaconis and Shahshahani prove that the
	eigenvalues for this random walk are the character ratios
	$\frac{\chi^{\lambda}(12)}{dim(\lambda)}$ each occurring with
	multiplicity $dim(\lambda)^2$. Hence Theorem \ref{ker} says
	that the spectrum of this random walk is asymptotically
	normal. Character ratios on transpositions also appear in work
	on the moduli space of curves \cite{EO}. In fact Kerov
	outlines a proof of the result. A full proof of the result
	appears in the marvelous paper \cite{IO}. Another approach is
	due to Hora \cite{H}, who exploited the fact that the $k$th
	moment of a Plancherel distributed character ratio is the
	chance that the random walk generated by random transpositions
	is at the identity after k steps (this follows from Lemma
	\ref{countsol} below). Both of these proofs establish asymptotic normality by the method of moments and use combinatorial methods to estimate the moments. Note also that there is no error term
	in Theorem \ref{ker}; in this paper one will be
	attained. Finally, we remark that Kerov (and then Hora) proves
	a much more general result-a multidimensional central limit
	theorem showing that character ratios evaluated on cycles of
	various lengths are asymptotically independent normal random
	variables.

	In this paper we prove the following result. Here $P(\cdot)$ denotes the probability of an event. 

\begin{theorem} \label{maintheorem} For $n \geq 2$ and all real $x_0$, \[ |P \left( \frac{n-1}{\sqrt{2}} \frac{\chi^{\lambda}(12)}{dim(\lambda)}  \leq x_0 \right) - \frac{1}{\sqrt{2 \pi}} \int_{-\infty}^{x_0} e^{-\frac{x^2}{2}} dx| \leq 40.1 n^{-1/4}. \]
\end{theorem}

	We conjecture that an upper bound of the form $C n^{-1/2}$
	holds (here $C$ is a constant). Even for simple problems,
	Stein's method can sometimes not (at least obviously) go
	beyond the $n^{-1/4}$ rate. The follow-up paper \cite{Fu2}
	gives a very interesting generalization of Theorem
	\ref{maintheorem} for Jack measure.

	Theorem \ref{maintheorem} will be a consequence of the
	following bound of Stein. Recall that if $W,W^*$ are random
	variables, they are called exchangeable if for all $w_1,w_2$,
	$P(W=w_1,W^*=w_2)$ is equal to $P(W=w_2,W^*=w_1)$. The
	notation $E^W(\cdot)$ means the expected value given $W$. Note from \cite{St1} that there are minor variations on Theorem \ref{steinbound} (and thus for Theorem \ref{maintheorem}) for $h(W)$ where $h$ is a bounded continuous function with bounded piecewise continuous derivative. For simplicity we only state the result when $h$ is the indicator function of an interval. 

\begin{theorem}

 \label{steinbound} (\cite{St1}) Let $(W,W^*)$ be an exchangeable pair
of real random variables such that $E^W(W^*) = (1-\tau)W$ with
$0<\tau<1$. Then for all real $x_0$, \begin{eqnarray*} & & |P(W \leq x_0) -
\frac{1}{\sqrt{2 \pi}} \int_{-\infty}^{x_0} e^{-\frac{x^2}{2}} dx|\\
& \leq & 2 \sqrt{E[1-\frac{1}{2 \tau}E^W (W^*-W)^2]^2} + (2 \pi)^{-
\frac{1}{4}} \sqrt{\frac{1}{\tau} E|W^*-W|^3}. \end{eqnarray*} \end{theorem}

	In order to apply Theorem \ref{steinbound} to study a
	statistic $W$, one clearly needs an exchangeable pair
	$(W,W^*)$ such that $E^W(W^*) = (1-\tau)W$ (this second
	condition can sometimes be weakened in using Stein's method
	\cite{RR}). Section \ref{harmonic} discusses the theory of
	harmonic functions on Bratelli diagrams and shows how it can
	be applied to generate a ``natural'' exchangeable pair
	$(W,W^*)$. The idea is to use a reversible Markov chain on the
	set of partitions of size $n$ whose stationary distribution is
	Plancherel measure, to let $\lambda^*$ be obtained from
	$\lambda$ by one step in the chain, and then set
	$(W,W^*)=(W(\lambda),W(\lambda^*))$. This construction also
	has the merit of being applicable to more general groups
	\cite{F},\cite{Fu2}.

	As we shall see in Section \ref{provebound}, we are quite
fortunate in that when $W=\frac{n-1}{\sqrt{2}}
\frac{\chi^{\lambda}(12)}{dim(\lambda)}$, it does happen that
$E^W(W^*) = (1-\tau)W$ for some $\tau$ (in fact
$\tau=\frac{2}{n+1}$). There are some other simplifications which
occur. For instance we will derive a simple upper bound for $2
\sqrt{E[1-\frac{1}{2 \tau}E^W (W^*-W)^2]^2}$, which is part of the
error term and consistent with then $C n^{-1/2}$ conjecture. To
appreciate the beauty of Stein's method, we note that there is only
one point in the proof of Theorem \ref{maintheorem} where we even use
an explicit formula for $W$, and this is in bounding $E|W^*-W|^3$. In
fact even this could be avoided by using the Cauchy Schwarz inequality
$E|Z|^3 \leq \sqrt{E(|Z|^2) E(|Z|^4)}$ where $Z$ is a random variable
(and this is useful for general finite groups), but the argument
presented involves longest increasing subsequences and we prefer
it. For results on other conjugacy classes and groups, see \cite{Fu3}.

	Section \ref{decompose} finds the eigenvalues and eigenvectors
	for the Markov chain underlying the construction of the
	exchangeable pair in Section \ref{harmonic}. This leads to a
	curious method for studying the decomposition of tensor
	products in the symmetric group. For example, let $V$ be the
	standard $n$-dimensional representation of the symmetric group
	$S_n$. We deduce that for $r$ sufficiently large (roughly
	$\frac{n^2 log(n)}{4}$), that for all $\lambda$, the
	multiplicity of the irreducible representation of type
	$\lambda$ of $S_n$ in the r-fold tensor product $V \otimes
	\cdots \otimes V$ is very close to $\frac{dim(\lambda)
	n^r}{n!}$. A follow-up paper \cite{F} stengthens this using
	card shuffling to show that $r$ roughly $n log(n)$ is
	sufficient and that $r$ of order $\frac{n log(n)}{2}$ is necessary.

	The precise organization of this paper is as follows. Section \ref{harmonic} uses the theory of harmonic functions to construct an exchangeable pair $(W,W^*)$. Section \ref{repreview} collects some lemmas we shall need from representation theory. We have done this for two reasons: first, to make the paper more readable by probabilists who work on Stein's method, and second, because there are a few new results and some non-standard facts. Section \ref{provebound} puts the pieces together and proves Theorem \ref{maintheorem}. Section \ref{decompose} finds the eigenvalues and eigenvectors for the Markov chain underlying the construction of the exchangeable pair in Section \ref{harmonic}, and applies it to obtain an asymptotic result about the decomposition of tensor products in the symmetric group. 

\section{Harmonic functions and exchangeable pairs} \label{harmonic} 

	To begin we recall the theory of harmonic functions on Bratelli diagrams. This is a beautiful subject with deep connections to probability theory and representation theory. Two excellent surveys are \cite{K2} and \cite{BO}.

	The basic set-up is as follows. One starts with a
Bratteli diagram; that is an oriented graded graph $\Gamma= \cup_{n \geq 0} \Gamma_n$ such that

\begin{enumerate}
\item $\Gamma_0$ is a single vertex $\emptyset$.
\item If the starting vertex of an edge is in $\Gamma_i$, then its end vertex is in $\Gamma_{i+1}$.
\item Every vertex has at least one outgoing edge.
\item All $\Gamma_i$ are finite.
\end{enumerate}

	For two vertices $\lambda, \Lambda \in \Gamma$, one writes $\lambda \nearrow \Lambda$ if there is an edge from $\lambda$ to
$\Lambda$. Part of the underlying data is a multiplicity function
$\kappa(\lambda,\Lambda)$. Letting the weight of a path in $\Gamma$ be
the product of the multiplicities of its edges, one defines the
dimension $dim(\Lambda)$ of a vertex $\Lambda$ to be the sum of the
weights over all maximal length paths from $\emptyset$ to $\Lambda$; $dim(\emptyset)$ is taken to be $1$. Given a Bratteli
diagram with a multiplicity function, one calls a function $\phi$ {\it
harmonic} if $\phi(\emptyset)=1$, $\phi(\lambda) \geq 0$ for all $\lambda \in
\Gamma$, and \[ \phi(\lambda) = \sum_{\Lambda: \lambda \nearrow
\Lambda} \kappa(\lambda,\Lambda) \phi(\Lambda).\] An equivalent
concept is that of coherent probability distributions. Namely a set
$\{M_n\}$ of probability distributions $M_n$ on $\Gamma_n$ is called
{\it coherent} if \[ M_{n}(\lambda) = \sum_{\Lambda: \lambda
\nearrow \Lambda} \frac{dim(\lambda)
\kappa(\lambda,\Lambda)}{dim(\Lambda)} M_{n+1}(\Lambda).\] The formula
allowing one to move between the definitions is $\phi(\lambda) =
\frac{M_n(\lambda)}{dim(\lambda)}$.

        One reason the set-up is interesting from the viewpoint of
probability theory is the fact that every harmonic function can be
written as a Poisson integral over the set of extreme harmonic
functions (i.e. the boundary). For the Pascal lattice (vertices of $\Gamma_n$ are pairs
$(k,n)$ with $k=0,1,\cdots,n$ and $(k,n)$ is connected to $(k,n+1)$
and $(k+1,n+1)$), this fact is the simplest instance of de Finetti's
theorem, which says that an infinite exchangeable sequence of $0-1$ random variables is a mixture of coin toss sequences for different
probabilities of heads. One of Kerov's insights in \cite{K2} is that one can prove Selberg type integral formulas by expressing interesting harmonic functions as integrals over the boundary. The paper \cite{BO} carries out this program in many cases.

	Before proceeding, let us indicate that Plancherel measure is
	a special case of this set-up. Here the lattice which one uses
	is the Young lattice, that is $\Gamma_n$ consists of all
	partitions of size $n$, and a partition of size $n$ is
	adjoined to a partition of size $n+1$ in $\Gamma_{n+1}$ if the
	partition of size $n+1$ can be obtained from the partition of
	size $n$ by adding a box. The multiplicity function
	$\kappa(\lambda,\Lambda)$ is equal to 1 on each edge. The
	dimension function $dim(\lambda)$ is simply
	$\frac{n!}{\prod_{x} h(x)}$, the dimension of the irreducible
	representation of the symmetric group parameterized by the
	partition $\lambda$. Then one can show that if $M_n$ is the
	Plancherel measure, the family of measures $\{M_n\}$ is
	coherent \cite{K2}. More generally it follows from Frobenius
	reciprocity that if $H$ is a subgroup of $G$, $M_H$ is
	Plancherel measure of $H$ and $M_G$ is Plancherel measure of
	$G$, and $\kappa(\lambda,\Lambda)$ is the multiplicity of
	$\lambda$ in the restriction of $\Lambda$ from $G$ to $H$,
	then \[ M_{H}(\lambda) = \sum_{\Lambda: \lambda \nearrow
	\Lambda} \frac{dim(\lambda)
	\kappa(\lambda,\Lambda)}{dim(\Lambda)} M_{G}(\Lambda).\]

	The relevance of these considerations is Proposition
	\ref{constructchain}. Parts 1 and 2 are implicit in
	\cite{K2}. Recall that a Markov chain $J$ on a finite set $X$
	is said to be reversible with respect to $\pi$ if $\pi(x)
	J(x,y) = \pi(y) J(y,x)$ for all $x,y$. It is easy to see that
	if $J$ is reversible with respect to $\pi$, then $\pi$ is
	stationary for $J$ (i.e. that $\pi(y) = \sum_{x \in X}
	\pi(x)J(x,y)$). Note that if $W$ is any statistic on $X$ and
	$W^*$ is obtained by evaluating the statistic after taking a
	step according to a Markov chain which is reversible with
	respect to $\pi$, then $(W,W^*)$ is an exchangeable pair under
	the probability measure $\pi$.

\begin{prop} \label{constructchain} Suppose that the family of measures $\{M_n\}$ is coherent for the Bratelli diagram. 
\begin{enumerate}
\item If $\lambda$ is chosen from the measure $M_n$, and one moves from $\lambda$ to $\Lambda$ with probability $\frac{dim(\lambda) M_{n+1}(\Lambda) \kappa(\lambda,\Lambda)}{dim(\Lambda) M_n(\lambda)}$, then $\Lambda$ is distributed according to the measure $M_{n+1}$.
\item If $\Lambda$ is chosen from the measure $M_{n+1}$, and one moves from $\Lambda$ to $\mu$ with probability $\frac{dim(\mu) \kappa(\mu,\Lambda)}{dim(\Lambda)}$, then $\mu$ is distributed according to the measure $M_n$.
\item The Markov chain $J$ on vertices in level $\Gamma_n$ of the Bratelli diagram given by moving from $\lambda$ to $\mu$ with probability \[ J(\lambda,\mu) = \frac{dim(\lambda) dim(\mu)}{M_n(\lambda)} \sum_{\Lambda \in \Gamma_{n+1}} \frac{M_{n+1}(\Lambda) \kappa(\lambda,\Lambda) \kappa(\mu,\Lambda)}{dim(\Lambda)^2}\] is reversible with stationary distribution $M_n$.
\item The Markov chain $J$ on vertices in level $\Gamma_n$ of the
Bratelli diagram given by moving from $\lambda$ to $\mu$ with
probability \[ J(\lambda,\mu) = \frac{M_n(\mu)}{dim(\lambda) dim(\mu)}
\sum_{\tau \in \Gamma_{n-1}} \frac{dim(\tau)^2 \kappa(\tau, \lambda)
\kappa(\tau,\mu)} { M_{n-1}(\tau)}\] is reversible with stationary
distribution $M_n$.
\end{enumerate}
\end{prop}

\begin{proof} For part 1, observe that \begin{eqnarray*} & & \sum_{\lambda \in \Gamma_n}
 M_n(\lambda) \frac{dim(\lambda) M_{n+1}(\Lambda)
 \kappa(\lambda,\Lambda)}{dim(\Lambda) M_n(\lambda)}\\ & = &
 \frac{M_{n+1}(\Lambda)}{dim(\Lambda)} \sum_{\lambda \in \Gamma_n} dim(\lambda) \kappa(\lambda,\Lambda)\\ & = &
 M_{n+1}(\Lambda).\end{eqnarray*} Note that the transition probabilities sum to 1 because the measures $\{M_n\}$ are coherent. Part 2 is similar and also uses the fact that $\{M_n\}$ is coherent. For part 3 reversibility is immediate from the definitions, provided that the transition probabilities for $J$ sum to 1. But $J$ is simply what one gets by moving up one level in the Bratelli diagram according to the transition
	mechanism of part 1, and then moving down according to the
	transition mechanism of part 2. Part 4 is similar; one first moves down the Bratelli diagram and then up. \end{proof}

	One can investigate more general Markov chains where one moves up or down by an amount $k$. This is done in Section \ref{decompose} where it is applied to tensor products. However in the application of Stein's method, we will only use the chain $J$ in part 3 of the proposition (the chain in part 4 would work as well). This chain simplifies in many cases of interest. We mention two of them. 

{\bf Example 1:} The first example is the Plancherel measure of $S_n$. We
already indicated how it fits in with harmonic functions on the Young
lattice. Let $parents(\lambda,\mu)$ denote the set of partitions above
both $\lambda,\mu$ in the Young lattice (this set has size $0$ or $1$
unless $\lambda=\mu$). Then
$J(\lambda,\mu)=\frac{dim(\mu)|parents(\lambda,\mu)|}{(n+1)dim(\lambda)}$. Lemma
\ref{parents} in the next section will use representation theory to
derive a more complicated (but useful) expression for $J(\lambda,\mu)$
as a sum over the symmetric group $S_n$.

{\bf Example 2:} The second example concerns cycles of random
permutations. As explained in \cite{K2}, the partitions of $n$ index
the conjugacy classes of the symmetric group and the probability
measure on partitions corresponding to the conjugacy class of a random
permutation is a coherent family of measures with respect to the
Kingman lattice. Here the underlying lattice is the same as the Young
lattice, but the multiplicity function $\kappa(\lambda,\Lambda)$ is
the number of rows of length $j$ in $\Lambda$, where $\lambda$ is
obtained from $\Lambda$ by removing a box from a row of length
$j$. Also $dim(\lambda)=\frac{n!}{\lambda_1 ! \cdots \lambda_l!}$
where $l$ is the number of rows of $\lambda$ and $\lambda_i$ is the
length of row $i$ of $\lambda$. So the Markov chain $J$ amounts to the
following. Letting $m_r(\lambda)$ denote the number of rows of length
$r$ in $\lambda$, first obtain $\Lambda$ by adding a box to a row of
length $r$ with chance $\frac{r m_r(\lambda)}{n+1}$, or to an empty
row with probability $\frac{1}{n+1}$. Then remove a box from a row of
$\Lambda$ of length $s$ with probability $\frac{s
m_s(\Lambda)}{n+1}$. The cycle structure of random permutations (and
more generally the Ewens sampling formula) is very well understood,
but still in future work we intend to revisit them by Stein's method
(using the exchangeable pairs in this section) and to study the
spectral properties of the Markov chains of Proposition
\ref{constructchain}.

\section{Irreducible characters of the symmetric group} \label{repreview}

	This section collects properties we will use about
	characters of irreducible representations of the symmetric
	group. For more background on this topic, see the book
	\cite{Sa}. Lemmas \ref{orth1}, \ref{mult}, and $\ref{real}$ are
	well known. Lemma \ref{countsol} is known but not well
	known. Lemmas \ref{induced}, \ref{parents} are elementary
	consequences of known facts but perhaps new.

	Throughout $\overline{z}$ denotes the complex conjugate of $z$, and $Irr(G)$ denotes the set of characters of irreducible representations of a finite group $G$. The notation $dim(\chi)$ denotes the dimension of the representation with character $\chi$.

\begin{lemma} \label{orth1} Let $C(g)$ be the conjugacy class of $G$ containing the element $g$. Then for $g \in G$, \[ \sum_{\chi \in Irr(G)} \chi(g) \overline{\chi(h)} \] is equal to $\frac{|G|}{|C(g)|}$ if $h,g$ are conjugate and is $0$ otherwise.
\end{lemma}

\begin{lemma} \label{mult} Let $\nu$ be an irreducible character of a finite group $G$, and $\chi$ any character of $G$. Then the multiplicity of $\nu$ in $\chi$ is equal to \[ \frac{1}{|G|} \sum_{g \in G} \nu(g) \overline{\chi(g)}. \] \end{lemma}

\begin{lemma} \label{real} The irreducible characters of the symmetric group $S_n$ are real valued. \end{lemma}

	Note that Lemma \ref{countsol} generalizes Lemma \ref{orth1}. 

\begin{lemma} (\cite{Sta}, Exercise 7.67) \label{countsol} Let $G$ be a finite group
 with conjugacy classes $C_1, \cdots, C_r$. Let $C_k$ be the conjugacy class of an element $w \in G$. Then the number of
 $m$-tuples $(g_1,\cdots,g_m) \in G^m$ such that $g_j \in C_{i_j}$ and $g_1 \cdots g_m=w$ is \[
 \frac{\prod_{j=1}^m |C_{i_j}|}{|G|} \sum_{\chi \in Irr(G)}
 \frac{1}{dim(\chi)^{m-1}} \chi(C_{i_1}) \cdots \chi(C_{i_m})
 \overline{\chi(C_k)} .\] \end{lemma}

	For the application of Stein's method we shall only need the case $k=1$ of Lemmas \ref{induced} and \ref{parents}.

\begin{lemma} \label{induced} Let $\chi$ be a character of the symmetric group $S_n$. Let $n_i(g)$ denote the number of cycles of $g$ of length $i$. Let $Res,Ind$ denote the operations of restriction and induction. Then for $k \geq 1$, \[ Res^{S_{n+k}}_{S_n}(Ind_{S_n}^{S_{n+k}}(\chi)) [g] = \chi(g) (n_1(g)+1) \cdots (n_1(g)+k).\]
\[ Ind_{S_{n-k}}^{S_n}(Res^{S_n}_{S_{n-k}}(\chi)) [g] = \chi(g) (n_1(g)) \cdots (n_1(g)-k+1).\]
\end{lemma}

\begin{proof} Let us prove the first assertion, the second being similar. It is well known (\cite{Sa}) that if $H$ is a subgroup of a finite group $G$, and $\chi$ is a character of $H$, then \[ Ind_H^G(\chi)[g] = \frac{1}{|H|} \sum_{t \in G \atop t^{-1}gt \in H} \chi(t^{-1}gt).\] Now apply this to $H=S_n$, $G=S_{n+k}$, using the fact that if two elements of $S_n$ are conjugate by an element of $S_{n+k}$, they are conjugate in $S_n$. Letting $Cent_G(g)$ denote the centralizer size of an element $g$ in a group $G$, it follows that the induced character at $g \in S_n$ is equal to \[ \chi(g) \frac{|Cent_{S_{n+k}}(g)|} {|Cent_{S_{n}}(g)|} = \chi(g) (n_1(g)+1) \cdots (n_1(g)+k).\] Here we have used the fact that $\prod_j j^{n_j} n_j!$ is the centralizer size in $S_{\sum_j jn_j}$ of an element with $n_j$ cycles of length $j$. \end{proof}

	In Lemma \ref{parents} we use the notation that $dim(\mu/\tau)$ is the number of paths in the Young lattice from $\tau$ to $\mu$ (or equivalently the number of ways of adding boxes one at a time to get from $\tau$ to $\mu$). We let $|\lambda|$ denote the size of a partition.

\begin{lemma} \label{parents} Let $\mu,\lambda$ be partitions of $n$. Let $n_1(g)$ denote the number of fixed points of a permutation $g$. Then \[ \sum_{\tau \atop |\tau|=n+k} dim(\tau / \lambda) dim(\tau / \mu) = \frac{1}{n!} \sum_{g \in S_n} \chi^{\mu}(g) \chi^{\lambda}(g) (n_1(g)+1) \cdots (n_1(g)+k)\] 
\[ \sum_{\tau \atop |\tau|=n-k} dim(\lambda / \tau) dim(\mu / \tau) = \frac{1}{n!} \sum_{g \in S_n} \chi^{\mu}(g) \chi^{\lambda}(g) (n_1(g)) \cdots (n_1(g)-k+1).\]
\end{lemma}

\begin{proof} We prove only the first part as the second part is similar. Both sides of this equation enumerate the multiplicity of $\mu$ in the representation $Res^{S_{n+k}}_{S_n}(Ind_{S_n}^{S_{n+k}}(\lambda))$. This multiplicity is equal to the left hand side by the branching rules for induction and restriction in the symmetric group (Theorem 2.8.3 in \cite {Sa}). That the right hand side computes the same multiplicity follows from Lemmas \ref{mult}, \ref{real}, and \ref{induced}. \end{proof}

	Note that as a special case of Lemma \ref{parents}, \[ |parents(\mu,\lambda)| = \frac{1}{n!} \sum_{g \in S_n} \chi^{\mu}(g) \chi^{\lambda}(g) (n_1(g)+1).\]

\section{Stein's method and Kerov's central limit theorem} \label{provebound}

	In this section we prove Theorem \ref{maintheorem}. Recall
	that $W(\lambda) = \frac{(n-1)}{\sqrt{2}}
	\frac{\chi^{\lambda}(12)}{dim(\lambda)}$ and that the
	Plancherel measure chooses a partition $\lambda$ with
	probability $\frac{dim(\lambda)^2}{n!}$. If $n=1$ we use the convention that $W=0$. Then Lemma
	\ref{countsol} implies the known fact that the mean and
	variance of $W$ are $0$ and $1-\frac{1}{n}$ respectively. This
	will also follow from Stein's method.

	In fact there is an explicit formula (due to Frobenius
\cite{Fr}) \[ \frac{\chi^{\lambda}(12)}{dim(\lambda)} = \frac{1}{{n
\choose 2}} \sum_i {\lambda_i \choose 2} - {\lambda_i' \choose 2} \]
where $\lambda_i$ is the length of row $i$ of $\lambda$ and
$\lambda_i'$ is the length of column $i$ of $\lambda$. We shall use
this formula only once.

	Throughout this section and the remainder of the paper, $W^*$ denotes $W(\lambda^*)$, where given $\lambda$, the partition $\lambda^*$ is $\mu$ with probability $J(\lambda,\mu)$ from Example 1 of Section \ref{harmonic}. Recall that $(W,W^*)$ is an exchangeable pair. We also use the notation that $\chi^{\lambda}$ denotes the character of the irreducible representation of the symmetric group parameterized by the partition $\lambda$. We let $|\lambda|$ denote the size of a partition.

	Proposition \ref{1strelate} shows that the hypothesis needed to apply the Stein method bound are satisfied. It also tells us that $W$ is an eigenvector for the Markov chain $J$, with eigenvalue $(1-\frac{2}{n+1})$. We shall generalize this observation in Section \ref{decompose}. 

\begin{prop} \label{1strelate} $E^W(W^*) = (1-\frac{2}{n+1}) W$. \end{prop}

\begin{proof} For $n=1$ we have that $W=0$ by convention. Otherwise, \begin{eqnarray*} E^{\lambda}(W^*) & = & \sum_{|\mu|= n} \frac{n-1}{\sqrt{2}} \frac{\chi^{\mu}(12)}{dim(\mu)} \frac{dim(\mu)}{(n+1) dim(\lambda)} |parents(\mu,\lambda)|\\
& = & \frac{n-1}{\sqrt{2}} \frac{1}{(n+1)dim(\lambda)} \sum_{|\mu| =n}
|parents(\mu,\lambda)| \chi^{\mu}(12). \end{eqnarray*} From the
representation theory of the symmetric group \cite{Sa},
$|parents(\mu,\lambda)|$ is equal to the multiplicity of $\chi^{\mu}$
in $Res^{S_{n+1}}_{S_n}(Ind_{S_n}^{S_{n+1}} (\chi^{\lambda}))$, since
inducing to $S_{n+1}$ corresponds to the possible ways of adding a box
to each corner of the partition $\lambda$, and restricting to $S_n$
corresponds to the possible ways of removing a corner box. Hence \[
E^{\lambda}(W^*) = \frac{n-1}{\sqrt{2}} \frac{1}{(n+1)dim(\lambda)}
Res^{S_{n+1}}_{S_n}(Ind_{S_n}^{S_{n+1}} (\lambda))[(12)] .\] Applying
Lemma \ref{induced}, this simplifies to \[ \frac{n-1}{\sqrt{2}}
\frac{\chi^{\lambda}(12)}{(n+1)dim(\lambda)} (n-1) =
(1-\frac{2}{n+1})W .\] Since $E^{\lambda}(W^*)$ depends on $\lambda$
only through $W$, it is equal to $E^W(W^*)$. \end{proof}

	As a consequence of Proposition \ref{1strelate}, we obtain a Stein's method proof that the mean $E(W)$ is equal to 0.

\begin{cor} $E(W)=0$. \end{cor} \begin{proof} Since the pair $(W,W^*)$ is exchangeable, $E(W^*-W)=0$. Using Proposition \ref{1strelate}, we see that \[ E(W^*-W)  =  E(E^W(W^*-W))
 = -\frac{2}{n+1} E(W) .\] Hence $E(W)=0$. \end{proof}

	Next we shall use Stein's method to compute $E^{\lambda}(W^*)^2$. Recall that this notation means the expected value of $(W^*)^2$ given $\lambda$. This will be useful for analyzing the error term in Theorem \ref{maintheorem}.  

\begin{prop} \label{forterm1} \begin{eqnarray*} E^{\lambda}((W^*)^2)
& = & (1-\frac{1}{n}) + \frac{2(n-1)(n-2)^2}{n(n+1)} \frac{\chi^{\lambda}((123))}{dim(\lambda)}\\
& &  +\frac{(n-1)(n-2)(n-3)^2}{2n(n+1)} \frac{\chi^{\lambda}((12)(34))}{dim(\lambda)}.\end{eqnarray*} Here we use the convention that if $n \leq 3$, $\chi^{\lambda}((12)(34))=0$, and if $n \leq 2$, $\chi^{\lambda}(123)=0.$ \end{prop} 

\begin{proof} If $n=1$ the result is clear. Otherwise, 
\begin{eqnarray*} E^{\lambda}(W^*)^2 & = & \frac{(n-1)^2}{2} \sum_{|\mu|=n} \frac{1}{n+1} |parents(\mu,\lambda)| \frac{dim(\mu)}{dim(\lambda)} \left(\frac{\chi^{\mu}(12)}{dim(\mu)} \right)^2\\
& = &  \frac{(n-1)^2}{2(n+1)} \frac{1}{dim(\lambda)} \sum_{|\mu|=n} |parents(\mu,\lambda)| \frac{\chi^{\mu}(12)^2}{dim(\mu)}. \end{eqnarray*} Applying Lemma \ref{parents}, this can be rewritten as \begin{eqnarray*}
& & \frac{(n-1)^2}{2(n+1)} \frac{1}{dim(\lambda)} \sum_{|\mu|=n} \frac{\chi^{\mu}(12)^2}{dim(\mu)} \frac{1}{n!} \sum_{g \in S_n} \chi^{\mu}(g) \chi^{\lambda}(g) (n_1(g)+1)\\
& = & \frac{(n-1)^2}{2(n+1)} \frac{1}{dim(\lambda)} \sum_{g \in S_n} \chi^{\lambda}(g) (n_1(g)+1) \frac{1}{n!} \sum_{|\mu|=n} \frac{\chi^{\mu}(12)^2 \chi^{\mu}(g)}{dim(\mu)}. \end{eqnarray*}

	The next step is to observe that using Lemmas \ref{real} and \ref{countsol}, one can compute the expression $\frac{1}{n!} \sum_{|\mu|=n} \frac{\chi^{\mu}(12)^2 \chi^{\mu}(g)}{dim(\mu)}$ for any permutation $g$. Indeed, it is simply $\frac{1}{{n \choose 2}^2}$ multiplied by the number of ordered pairs $(\tau_1,\tau_2)$ of transpositions whose product is $g$. Thus when this expression is non-zero, there are 3 possibilities for the cycle type of $g$: the identity, a 3-cycle, or a product of two 2-cycles on disjoint symbols. In all cases it is elementary to enumerate the number of pairs $(\tau_1,\tau_2)$, and these 3 possibilities yield the 3 terms in the statement of the proposition. 
\end{proof}

	As a consequence of Lemma \ref{forterm1}, we compute $Var(W)$. 

\begin{cor} $Var(W)=(1-\frac{1}{n})$. \end{cor}

\begin{proof} Since $W$ has mean $0$, $Var(W)=E(W^2)$. Since $W$ and $W^*$ have the same distribution, it follows that \[ E(W^2)=E(W^*)^2=E(E^{\lambda}(W^*)^2).\] The quantity $E^{\lambda}(W^*)^2$ was computed in Proposition \ref{forterm1} as a sum of three terms. The Plancherel measure average of the first term is $(1-\frac{1}{n})$. The Plancherel measure averages of the other terms both vanish. Indeed, if $g$ is any nonidentity element of the symmetric group, the Plancherel average of the function $\frac{\chi^{\lambda}(g)}{dim(\lambda)}$ is equal to \[ \frac{1}{n!} \sum_{\lambda} dim(\lambda) \chi^{\lambda}(g)
= \frac{1}{n!} \sum_{\lambda} \chi^{\lambda}(1) \chi^{\lambda}(g), \]
which vanishes by Lemma \ref{orth1}. \end{proof}

	In order to prove Theorem \ref{maintheorem}, we have to
	analyze the error terms in Theorem \ref{steinbound}. To begin
	we study \[ E\left(-1+\frac{n+1}{4} E^{\lambda}(W^*-W)^2
	\right)^2, \] obtaining an exact formula. From a well known
	general principle used in \cite{St1} (see Lemma 5 of
	\cite{Fu4} for a proof) the fact that $W$ is determined by
	$\lambda$ implies that \[ E[E^W(W^*-W)^2]^2 \leq
	E[E^{\lambda}(W^*-W)^2]^2 .\] Hence Proposition \ref{term1}
	gives an upper bound on \[ E\left(-1+\frac{n+1}{4}
	E^W(W^*-W)^2 \right)^2. \]

\begin{prop} \label{term1} \[ E \left(-1+\frac{n+1}{4} E^{\lambda}(W^*-W)^2 \right)^2 = \frac{3n^2-5n+6}{4n^3}.\] \end{prop}

\begin{proof} First observe (using Proposition \ref{1strelate} in the second equality)
 that \[ E^{\lambda}(W^*-W)^2 = W^2 - 2WE^{\lambda} W^* + E^{\lambda}(W^*)^2 =
 (\frac{4}{n+1}-1)W^2 + E^{\lambda}(W^*)^2.\] Combining this with Proposition
 \ref{forterm1}, it follows that $-1+\frac{n+1}{4} E^{\lambda}(W^*-W)^2$ is
 equal to $A+B+C+D$ where \begin{enumerate}
\item $A= \left( \frac{n+1}{4}(1-\frac{1}{n}) - 1 \right)$
\item $B= \frac{(n-1)(n-2)(n-3)^2}{8n} \frac{\chi^{\lambda}((12)(34))}{dim(\lambda)}$
\item $C= \frac{(n-1)(n-2)^2}{2n} \frac{\chi^{\lambda}(123)}{dim(\lambda)}$
\item $D= -\frac{n+1}{4}(1-\frac{4}{n+1}) \frac{(n-1)^2}{2} \left( \frac{\chi^{\lambda}(12)}{dim(\lambda)} \right)^2$
\end{enumerate} 

	We need to compute the Plancherel measure average of
$(A+B+C+D)^2$. Since $A^2$ is constant, the average of $A^2$ is
$\left( \frac{n+1}{4}(1-\frac{1}{n}) - 1 \right)^2$. The Plancherel
averages of $B^2$,$C^2$ can both be computed using Lemma
\ref{orth1}. One gets $\frac{(n-1)(n-2)(n-3)^3}{8n^3}$ and
$\frac{3(n-1)(n-2)^3}{4n^3}$ respectively. To compute the Plancherel
average of $D^2$ one uses Lemma \ref{countsol} to reduce the
computation to counting the number of ordered triples
$(\tau_1,\tau_2,\tau_3)$ of transpositions whose product is a
transposition, which is easy to do. The Plancherel averages of $2AB$,
$2AC$, $2BC$ are all 0 by Lemma \ref{orth1}. The Plancherel average of
$2AD$ is computed by Lemma \ref{orth1}. The Plancherel average of
$2BD$ is reduced to counting the number of ordered pairs of
transpositions whose product consists of two disjoint 2 cycles by
Lemma \ref{countsol}, and similarly the average of $2CD$ is reduced to
counting the number of ordered pairs of transpositions whose product
is a 3 cycle. Thus all of the enumerations are elementary and adding
up the terms yields the proposition. \end{proof}

	The final ingredient needed to prove Theorem \ref{maintheorem}
	is an upper bound on $E|W^*-W|^3$. Typically this is the
	crudest term in applications of Stein's method. It is the only
	point in the argument where we use Frobenius' explicit formula
	for $W$ (and as noted in the introduction even this could be
	avoided).

\begin{prop} \label{term2} \[ E|W^*-W|^3 \leq \left(\frac{4e \sqrt{2}}{\sqrt{n}} \right)^3 + 2e^{-2e \sqrt{n}} (2 \sqrt{2})^3 \]
\end{prop} 

\begin{proof} From Frobenius' formula, \[ W= \frac{\sqrt{2}}{n} \sum_i {\lambda_i \choose 2} - {\lambda_i' \choose 2} .\] Given this and the way that $\lambda^*$ is constructed from $\lambda$, it follows that \[|W^*-W| \leq \frac{\sqrt{2}}{n} 2 max(\lambda_1,\lambda_1').\] Indeed, suppose that $\lambda^*$ is obtained from $\lambda$ by moving a box from row $a$ and column $b$ to row $c$ and column $d$. Then \[ W^*-W = \frac{\sqrt{2}}{n} \left( \lambda_c + \lambda_b' - \lambda_a - \lambda_d' \right).\] 

	Suppose that $\lambda_1$ (the size of the first row of $\lambda$) and $\lambda_1'$ (the size of the first column of $\lambda$) are both at most $2e \sqrt{n}$. Then by the previous paragraph \[ |W^*-W| \leq \frac{4e \sqrt{2}}{\sqrt{n}}.\] 

	Next note by the first paragraph, even if $\lambda_1>2e \sqrt{n}$ or $\lambda_1'>2e \sqrt{n}$ occurs, $|W^*-W| \leq 2 \sqrt{2}$. We claim that the chance that at least one of the events $\lambda_1>2e \sqrt{n}$ or $\lambda_1'>2e \sqrt{n}$ occurs is at most $2e^{-2e \sqrt{n}}$. Indeed, it is a simple lemma proved on page 7 of \cite{Ste} that the chance that the longest increasing subsequence of a random permutation is at least $2e \sqrt{n}$ is at most $e^{-2e \sqrt{n}}$. But it follows from the Robinson-Schensted-Knuth correspondence (see for instance \cite{Sa}) that when $\lambda$ is chosen from the Plancherel measure of the symmetric group on $n$ symbols, both $\lambda_1$ and $\lambda_1'$ have the same distribution as the longest increasing subsequence of a random permutation on $n$ symbols. So as claimed, the chance that at least one of the events $\lambda_1>2e \sqrt{n}$ or $\lambda_1'>2e \sqrt{n}$ occurs is at most $2e^{-2e \sqrt{n}}$. 

	Combining these observations proves the proposition. \end{proof}

	To close this section, we prove Theorem \ref{maintheorem}.

\begin{proof} We use Theorem \ref{steinbound}, which is applicable with $\tau=\frac{2}{n+1}$ by Propositions \ref{constructchain} and \ref{1strelate}. Since $n \geq 2$, Proposition \ref{term1} implies that the term \[ 2 \sqrt{E[1-\frac{1}{2 \tau}E^W
(W^*-W)^2]^2} \] in Stein's bound is at most $\sqrt{3} n^{-1/2} \leq \sqrt{3} n^{-1/4}$. So using Proposition \ref{term2} and the fact that $e^{-2e \sqrt{n}} \leq n^{-3/2}$, the bound of Theorem \ref{steinbound} becomes \[ \sqrt{3} n^{-1/4} + (2 \pi)^{-1/4} \sqrt{\frac{n+1}{2} \left( (\frac{4e\sqrt{2}}{\sqrt{n}})^3 + \frac{2 (2 \sqrt{2})^3}{n^{3/2}} \right)} \leq 40.1 n^{-1/4}.\] \end{proof}

\section{Asymptotic multiplicities in tensor products of representations} \label{decompose}

	This section gives an intruiging method for understanding some
	asymptotic aspects of the decomposition of tensor products of
	certain representations of the symmetric group. A recent paper
	which investigates this topic from the viewpoint of free
	probability theory is \cite{Bi}. However the methods and
	results are completely different from those presented
	here. The main new idea of this section is to use spectral
	theory of Markov chains. This is further developed in
	\cite{F} (connections with card shuffling, generalization to
	other groups) and in \cite{Fu2} (generalizations to arbitrary
	representations of finite groups).

	Throughout this section $X$ is the set of partitions of size $n$, endowed with Plancherel measure $\pi$. We consider the space of real valued functions $\ell^2(\pi)$ with the norm \[ ||f||_2 = \left( \sum_x |f(x)|^2 \pi(x) \right)^{1/2}.\] If $J(x,y)$ is the transition rule for a Markov chain, the associated Markov operator (also denoted by $J$) on $\ell^2(\pi)$ is given by $Jf(x) = \sum_y J(x,y) f(y)$. Let $J^r(x,y)=J_x^r(y)$ denote the chance that the Markov chain started at $x$ is at $y$ after $r$ steps.

	If the Markov chain with transition rule $J(x,y)$ is reversible with respect to $\pi$ (i.e. $\pi(x) J(x,y) = \pi(y) J(y,x)$ for all $x,y$), then the operator $J$ is self adjoint with real eigenvalues \[ -1 \leq \beta_{min}=\beta_{|X|-1} \leq \cdots \leq \beta_1 \leq \beta_0=1.\] Let $\psi_i$ ($i=0,\cdots,|X|-1$) be an orthonormal basis of eigenfunctions such that $J\psi_i = \beta_i \psi_i$ and $\psi_0 \equiv 1$. Define $\beta = max \{|\beta_{min}|,\beta_1\}$.

	The total variation distance between the measures $J_x^r$ and $\pi$ is defined by $||J_x^r- \pi||_{TV} = \frac{1}{2} \sum_{y} |J_x^r(y)-\pi(y)|$. (The application to representation theory does not require total variation distance but this concept is useful for understanding the Markov chain $J$; see Theorem \ref{main2}).

	The following lemma is well-known; for a proof see \cite{DSa}. Part 1 is essentially Jensen's inequality. 

\begin{lemma} \label{classic}
\begin{enumerate}
\item $2 ||J_x^r - \pi||_{TV} \leq ||\frac{J_x^r}{\pi}-1||_2$.
\item $J^r(x,y) = \sum_{i=0}^{|X|-1} \beta_i^r \psi_i(x) \psi_i(y) \pi(y)$.
\item $||\frac{J_x^r}{\pi}-1||_2^2 = \sum_{i=1}^{|X|-1} \beta_i^{2r} |\psi_i(x)|^2 \leq \frac{1-\pi(x)}{\pi(x)} \beta^{2r}$.
\end{enumerate}
\end{lemma}

	Next we specify a class of Markov chains to which we will apply Lemma \ref{classic}. Recall from Proposition \ref{constructchain} that there is a natural transition mechanism for moving up one step in the Young lattice or down one step in the Young lattice. For $1 \leq k \leq n$, let $J(k)$ be the Markov chain which from a partition of $n$ first moves down $k$ steps in the Young lattice (one at a time according to part 2 of Proposition \ref{constructchain}), and then moves back up $k$ steps in the Young lattice (one at a time according to part 1 of Proposition \ref{constructchain}). Proposition \ref{diagonalize} finds the eigenvalues and eigenfunctions for the corresponding operators $J(k)$. The word Mult. stands for multiplicity.

\begin{prop} \label{diagonalize} The eigenvalues and eigenfunctions of the operator $J(k)$ on partitions of size $n$ are indexed by conjugacy classes $C$ of the symmetric group on $n$ symbols.
\begin{enumerate}
\item Letting $n_1(C)$ denote the number of fixed points of the class $C$, the eigenvalue parameterized by $C$ is $\frac{(n_1(C))(n_1(C)-1) \cdots (n_1(C)-k+1)}{n(n-1) \cdots (n-k+1)}$.
\item The orthonormal basis of eigenfunctions $\psi_C$ is defined by $\psi_C(\lambda)=|C|^{\frac{1}{2}} \frac{\chi^{\lambda}(C)}{dim(\lambda)}$.
\end{enumerate}
\end{prop}

\begin{proof} From Proposition \ref{constructchain}, the chance that the Markov
chain $J(k)$ moves from $\lambda$ to $\mu$ is easily seen to be \[ \frac{dim(\mu)}{ (n)(n-1) \cdots (n-k+1) dim(\lambda)} \sum_{\tau \atop |\tau|=n-k} dim(\lambda/\tau) dim(\mu/\tau).\] From the branching rules of irreducible representations of the symmetric group \cite{Sa}, this is   \[ \frac{dim(\mu)}{dim(\lambda)} \frac{Mult. \ \mu \
in \ Ind^{S_n}_{S_{n-k}} Res^{S_n}_{S_{n-k}} (\lambda)}{n(n-1)\cdots
(n-k+1)}.\] Now observe that the $\psi_C(\lambda)=|C|^{\frac{1}{2}}
\frac{\chi^{\lambda}(C)}{dim(\lambda)}$ is an eigenfunction with the asserted eigenvalue because
\begin{eqnarray*} & & \frac{|C|^{\frac{1}{2}}}{(n) \cdots (n-k+1) dim(\lambda)} \sum_{\mu} \chi^{\mu}(C) \cdot  Mult. \ \mu \ in \ Ind^{S_n}_{S_{n-k}}
Res^{S_n}_{S_{n-k}}(\lambda) \\ & = & \frac{|C|^{\frac{1}{2}}}{n (n-1) \cdots (n-k+1)
dim(\lambda)} Ind^{S_n}_{S_{n-k}} Res^{S_n}_{S_{n-k}}(\lambda) [C]\\ & = &
\frac{(n_1(C))(n_1(C)-1) \cdots (n_1(C)-k+1)}{n(n-1)\cdots (n-k+1)}
\frac{|C|^{\frac{1}{2}} \chi^{\lambda}(C)}{dim(\lambda)}. \end{eqnarray*} The last equality is Lemma \ref{induced}. The fact that
$\psi_C$ are orthonormal follows from Lemma \ref{orth1}. Being
linearly independent, they are a basis for $\ell ^2(\pi)$ since the
number of conjugacy classes of $S_n$ is equal to the number of partitions of
$n$. \end{proof}

Remarks:
\begin{enumerate}
\item A similar argument shows that the chain which moves up $k$ steps (one at a time) and then down $k$ steps (one at a time) has the $\psi_C$ as an orthonormal basis with eigenvalues $\frac{(n_1(C)+1)(n_1(C)+2) \cdots (n_1(C)+k)}{(n+1)(n+2) \cdots (n+k)}.$ The results in the remainder of this section could be applied to that Markov chain as well.

\item It is remarkable that the set of eigenvalues of the Markov chain
$J(1)$ is precisely the set of eigenvalues of the top to random
shuffle (see \cite{DFP} for background on top to random
shuffles). This is not coincidence; the paper \cite{F} explains this
and more general facts in terms of the Robinson-Schensted-Knuth
correspondence.
\end{enumerate}

	Theorem \ref{main2} studies the convergence properties of the Markov chains $J(k)$. Throughout all logs are taken base e, as usual.
	
\begin{theorem} \label{main2} Let $(n)$ denote
 the partition which consists of one part of size $n$ (corresponding
 to the trivial representation of $S_n$). Let $\beta=\frac{(n-k)(n-k-1)}{(n)(n-1)}$. Then for $n \geq 3, 1 \leq k < n$, \[
 2||J(k)_{(n)}^r-\pi||_{TV} \leq \left( \sum_{\lambda \atop
 |\lambda|=n} |\frac{J(k)^r((n),\lambda)}{\pi(\lambda)}-1|^2
 \pi(\lambda) \right)^{1/2} \leq \sqrt{n!} \left( \beta
 \right)^{r}.\] Thus for $r > \frac{n log(n) +
 2c}{2log(\frac{1}{\beta})}$, \[ ||J(k)_{(n)}^r-\pi||_{TV} \leq
 \frac{(2 \pi)^{\frac{1}{4}}}{2} e^{-c}.\]
\end{theorem}

\begin{proof} From parts 1 and 3 of Lemma \ref{classic}, to prove the first assertion it is enough to show that for the chain $J(k)$, $\beta=\frac{(n-k)(n-k-1)}{(n)(n-1)}$ and that $\frac{1-\pi((n))}{\pi((n))} \leq n!$. The value of $\beta$ follows from Proposition \ref{diagonalize} (note that a permutation on n symbols cannot have exactly n-1 fixed points) and the inequality holds since $\pi((n))=\frac{1}{n!}$. This proves the first assertion. 

	The second assertion follows from the Stirling formula bound
\cite{Fe} \[ n! \leq \sqrt{2 \pi}
e^{-n+\frac{1}{12n}+(n+\frac{1}{2})log(n)}, \] for then
\begin{eqnarray*} \sqrt{n!} \left( \beta \right)^{r} &
\leq & (2 \pi)^{\frac{1}{4}} e^{r
log(\beta)-\frac{n}{2}+\frac{1}{24n}+\frac{1}{2}(n+\frac{1}{2})log(n)}\\
& \leq & (2 \pi)^{\frac{1}{4}} e^{r log(\beta)+\frac{n
log(n)}{2}}. \end{eqnarray*} \end{proof}

	Now we give the application to representation theory. We remind the reader that $Ind_{S_{n-1}}^{S_n}(1)$ is the defining representation of the symmetric group (i.e. the n dimensional permutation representation on the symbols $1,\cdots,n$), and hence corresponds to the case $k=1$ below. As usual, $\pi(\lambda)=\frac{dim(\lambda)^2}{n!}$ denotes the Plancherel measure of the symmetric group, and the word Mult. is short for multiplicity. We let $\otimes^r$ denote the operation of taking the r-fold tensor product of a representation of a symmetric group $S_n$ (yielding the sum of various irreducible representations of $S_n$).

\begin{theorem} \label{main3} Let $\beta=\frac{(n-k)(n-k-1)}{(n)(n-1)}$. Then for $n \geq 3, 1 \leq k < n$, \[  \sum_{\lambda
 \atop |\lambda|=n} |\frac{n! [Mult. \ \lambda \ in \ \otimes^r
 Ind_{S_{n-k}}^{S_n} (1)]}{dim(\lambda) \left( (n)\cdots (n-k+1)
 \right)^r}-1|^2 \pi(\lambda) \leq n! \left(\beta \right)^{2r}.\] For
 $r> \frac{n log(n) + c}{2log(\frac{1}{\beta})}$, this is at most
 $\sqrt{2 \pi} e^{-c}$. \end{theorem}

\begin{proof}
 Let $(n)$ be the partition which consists of one row of size $n$, and
 let $n_1(C)$ denote the number of fixed points of a conjugacy class
 $C$. We now consider the quantity $||\frac{J(k)_{(n)}^r}{\pi}-1||_2^2$. One on
 hand, by part 2 of Lemma \ref{classic}, we know that \begin{eqnarray*} & & J(k)_{(n)}^r
 (\lambda)\\ & = & dim(\lambda) \sum_{C} \left( \frac{(n_1(C)) \cdots
 (n_1(C)-k+1)}{n \cdots (n-k+1)} \right) ^r \frac{|C| \chi^{\lambda}(C)}{n!}\\
& = & \frac{dim(\lambda)}{\left( n \cdots (n-k+1) \right) ^r} \frac{1}{n!} \sum_{g \in S_n} \left( n_1(g) \cdots
 (n_1(g)-k+1) \right)^r \chi^{\lambda}(g). \end{eqnarray*} By Lemmas \ref{mult} and \ref{induced}, and the
 fact that the character of a tensor product is the product of the
 characters, this is precisely \[ \frac{dim(\lambda)}{\left( n(n-1) \cdots
 (n-k+1) \right) ^r} [Mult. \ \lambda \ in \ \otimes^r
 Ind_{S_{n-k}}^{S_n} (1)].\] Thus \[ ||\frac{J(k)_{(n)}^r}{\pi}-1||_2^2 = 
 \sum_{\lambda \atop |\lambda|=n} |\frac{ n! [Mult. \ 
 \lambda \ in \ \otimes^r Ind_{S_{n-k}}^{S_n} (1)]}{dim(\lambda) \left((n)\cdots(
 n-k+1) \right)^r}-1|^2 \pi(\lambda).\] The result now follows
 by the upper bound on $||\frac{J(k)_{(n)}^r}{\pi}-1||_2^2$ in Theorem \ref{main2}. \end{proof}

\section{Acknowledgements} The author was partially supported by National Security Agency grant MDA904-03-1-0049.

\end{document}